\documentclass[11pt,reqno]{amsart}
\usepackage{fonttable}
\usepackage{graphicx,color}
\usepackage{amssymb,amscd,latexsym, mathrsfs, epsfig}
\usepackage[all]{xy}
\usepackage[latin1]{inputenc}
\usepackage{mathdots}
\newtheorem{thm}{Theorem}[section]
\newtheorem{lem}[thm]{Lemma}
\newtheorem{pro}[thm]{Proposition}
\newtheorem{rem}[thm]{Remark}

\newtheorem{cor}[thm]{Corollary}
\newtheorem{defi}[thm]{Definition}

\newcommand{\bp}{{\bf P}}

\newcommand{\filt}[6]{
\[
\begin{xy}
\xymatrix@R20pt@C20pt{
&\mathbb{C}^3&\\\langle #1,#2 \rangle\ar[ru]&\langle #3,#4\rangle\ar[u]&\langle #5,#6\rangle\ar[lu]\\\langle #1\rangle\ar[u]&\langle #3\rangle\ar[u]&\langle #5\rangle\ar[u]}
\end{xy}
\]
}
\newcommand{\quiverext}[6]{\[
\begin{xy}
\xymatrix@R20pt@C20pt{
(#1)&&(#2)\ar[ld]^{#5}\ar[ll]^{#4}\\&(#3)\ar[lu]^{#6}}
\end{xy}
\]
}
\newcommand{\quiverhom}[6]{\[
\begin{xy}
\xymatrix@R20pt@C20pt{
(#1)&&(#2)\ar@{--}[ld]^{#5}\ar[ll]^{#4}\\&(#3)\ar[lu]^{#6}}
\end{xy}
\]
}

\newcommand\Zn{\mathbb{Z}}
\newcommand\Nn{\mathbb{N}}
\newcommand\Cn{\mathbb{C}}

\newcommand\cD{\mathcal D}
\newcommand\vm{v_{\max}}
\newcommand\es{\emptyset}
\newcommand\dv{\mathrm{div}}

\newcommand{\Hom}{\mathrm{Hom}}

\def\bd{\mathbf{ d}}
\def\st{\mathrm{st}}
\DeclareSymbolFont{symbolsC}{U}{pxsyc}{m}{n}
\DeclareMathSymbol{\Perp}{\mathrel}{symbolsC}{121}

\begin{document}
\title{Moduli spaces of point configurations and plane curve counts}

\begin{abstract} We identify certain Gromov-Witten invariants counting rational curves with given incidence and tangency conditions with the Betti numbers of moduli spaces of point configurations in projective spaces. On the Gromov-Witten side, S. Fomin and G. Mikhalkin established a recurrence relation via tropicalization, which is realized on the moduli space side using Donaldson-Thomas invariants of subspace quivers.\end{abstract}

\author{Markus Reineke}
\address{Markus Reineke: Ruhr-Universit\"at Bochum, Universit\"atsstra{\ss}e 150, 44780 Bochum, Germany}
\email{markus.reineke@ruhr-uni-bochum.de}
\author{Thorsten Weist}
\address{Thorsten Weist: Bergische Universit\"at Wuppertal, Gau\ss str.\ 20, 42097 Wuppertal, Germany}
\email{weist@uni-wuppertal.de}

\maketitle

\section{Introduction}

In this paper, we will point out a numerical relation between to seemingly unrelated geometries. On the one hand, we consider certain Gromov-Witten invariants, more precisely counts of irreducible rational curves in the projective plane satisfying suitable incidence and tangency conditions. The computation of such invariants is a priori of intersection theoretic nature. On the other hand, we consider moduli spaces of ordered point configurations in projective spaces up to projective symmetries, such spaces being classical examples of Geometric Invariant Theory. The invariants of these spaces we consider are of purely topological nature, namely, we are interested in their Euler characteristic, or, more generally, their Poincar\'e polynomials.\\[1ex]
Our main result, Corollary \ref{mainresult} below, is an equality between these invariants. Despite its rather elementary formulation, no direct geometric proof is available at the moment. Instead, we derive this equality by computation of both sides from the same recursive formula, see Theorem \ref{r2}.\\[1ex]
On the Gromov-Witten side, such a recursion is provided by work of S.~Fomin and G.~Mikhalkin which uses tropicalization and explicit combinatorial enumeration. On the moduli space side, we use an interpretation of the relevant geometries as moduli spaces of quiver representations, and derive the recursion from a DT/PT-type equation relating Poincare\'e polynomials of framed moduli spaces and motivic Donaldson-Thomas invariants of quivers with stability, together with a duality property and relations between framed and unframed quiver moduli spaces which are special to the present situation.\\[1ex]
At first view the relation established here is similar to the GW/Kronecker correspondence \cite{rw} which also provides a relation between Gromov-Witten invariants of toric surfaces and topological invariants of moduli spaces of quiver representations; however, there are substantial differences in the numerical parameters involved preventing this previous work to be applicable to the present relation. On a more technical level, whereas \cite{rw} relies on the wall-crossing formula for topological invariants of quiver moduli (and avoiding Donaldson-Thomas invariants), the present work makes substantial use of properties of Donaldson-Thomas invariants for quivers developed in recent years.\\[1ex]
Although the methods which will be developed here are rather specially taylored towards the claimed equality of invariants, they nevertheless provide some further applications, namely a new derivation of a formula for the Poincar\'e polynomials of moduli spaces of point configurations in the projective line, a rather curious identity between Gromov-Witten invariants, and a conjectural recursive formula for Block-G\"ottsche invariants generalizing the recursion of \cite{fm}.\\[1ex]
We also discuss the similarities and differences between the combinatorial objects which can be extracted from the two geometries, namely labeled floor diagrams on the Gromov-Witten side and stable trees on the moduli space side; it is remarkable that the numerical correspondence, although it can be formulated on a purely combinatorial level, apparently cannot be derived by combinatorial means only.\\[1ex]
We hope that the present work can serve as a starting point for establishing more general correspondences between Gromov-Witten geometries on the one hand and quiver moduli on the other hand, for example by suitably interpreting more of the invariants of \cite{fm} in this way, or by unifying the current picture with the setup of \cite{rw}. However, the identification of the relevant quiver moduli and the mechanisms leading to such numerical correspondences are missing at the moment.\\[1ex]
In the rest of this section, we will first give precise definitions of the relevant geometries, and formulate the numerical relations between invariants. In Section \ref{quiver}, we will recall in detail all relevant methods on quiver moduli. These will be applied,  specialized and complemented in Section \ref{point}. This allows us to give short proofs of the main results in Section \ref{app}. The correspondence of combinatorial objects will be described and discussed in Section \ref{comb}.

\subsection{Moduli spaces of point configuration}

We first recall the definition of moduli spaces of ordered point configurations in projective spaces up to projective symmetries, which is one of the classical examples of Geometric Invariant Theory \cite[Chapter 3]{MFK}.\\[1ex]
We fix nonnegative integers $d$ and $m$. An $m$-tuple $(p_1,\ldots,p_m)\in(\mathbb{P}^{d-1})^m$ of points in $\mathbb{P}^{d-1}$ is called semistable 
if, for every non-empty proper subset $I\subset[m]=\{1,\ldots,m\}$, the dimension of the linear subspace $U_I\subset\mathbb{P}^{d-1}$ generated by the $p_i$ for $i\in I$ fulfills
$$\dim U_I\geq \frac{d}{m}|I|-1
.$$

Denote by $(\mathbb{P}^{d-1})^m_{\mathrm{(s)st}}$ the (semi-)stable locus in the space $(\mathbb{P}^{d-1})^m$ of all $m$-tuples of points in $\mathbb{P}^{d-1}$, on which the projective linear group ${\rm PGL}_d$ acts diagonally.  We define
$$M_{d,m}=(\mathbb{P}^{d-1})^m_{\mathrm{sst}}//\mathrm{PGL}_d$$
as the quotient parametrizing closed ${\rm PGL}_d$-orbits, and call it the moduli space of semistable point configurations.\\[1ex]
The variety $M_{d,m}$ is non-empty if and only if $d<m$, in which case it is irreducible and projective of dimension $(m-d-1)(d-1)$. If $d$ and $m$ are coprime, $M_{d,m}$ is smooth, otherwise it is typically singular at properly semistable configurations.\\[1ex]
Only few of these moduli spaces identify with classical spaces of projective geometry: $M_{2,4}$ identifies with the projective line via the cross-ratio of four points, $M_{2,5}$ is a del Pezzo surface of degree five, and $M_{2,6}$ is isomorphic to the Segre cubic threefold. Explicit coordinates in the cases $d=2$ or $m=6$ are derived in \cite{HMSV,HMSV1}. An efficient formula for the Betti numbers of the spaces $M_{d,m}$ in the case of coprime $d$ and $m$ is given in \cite[Section 7]{HNS}. Recently, the spaces $M_{2,2k}$ were shown to be examples of spaces admitting two non-isomorphic small resolutions \cite{FR}.

\subsection{Gromov-Witten invariants}

Following \cite{fm}, we denote by $N_{d,g}$ the number of irreducible curves in the complex projective plane $\mathbb P^2$ of degree $d$ and genus $g$ passing through $3d+g-1$ fixed points in general position. For $g=0$, these numbers are determined recursively by the famous Kontsevich recursion \cite[5.17]{KM}:
$$N_{d,0}=\sum_{k+l=d}N_{k,0}N_{l,0}k^2l(l{3d-4\choose 3k-2}-k{3d-4\choose 3k-1}).$$ 
For fixed partitions $\lambda$ and $\rho$, the relative Gromov-Witten invariant $N_{d,g}(\lambda,\rho)$ is defined as the number of irreducible curves of  degree $d$ and genus $g$ passing through $2d-1+g+l(\rho)$ fixed points in general position which satisfy tangency conditions, with respect to a given line $L$, which are described by the partitions $\rho$ and $\lambda$, see \cite[Definition 3.15]{fm} for a precise definition. Note that the points corresponding to the partition $\lambda$ are fixed whence those corresponding to $\rho$ may vary. We are mostly interested in the following two special cases:
\begin{itemize}
\item $N_{d,0}((d),\emptyset)$ which counts irreducible rational degree $d$ curves in $\mathbb{P}^2$ passing through a given generic configuration of $2d-1$ points and having order $d$ tangency to a given line at a given point,
\item $N_{d,0}(\emptyset,(d))$ which counts irreducible rational degree $d$ curves in $\mathbb{P}^2$ passing through a given generic configuration of $2d$ points and having order $d$ tangency to a given line at an unspecified point.
\end{itemize}

\subsection{Main results}\label{results}





For a connected complex variety $X$, define

$$P_X(q)=(-q^\frac{1}{2})^{-\dim X}\sum_i\dim {\rm IH}^i(X,\mathbb{Q})(-q^\frac{1}{2})^i\in\mathbb{Z}[q^{\pm\frac{1}{2}}]$$
as a shifted version of its Poincar\'e polynomial in intersection cohomology. For $X$ a projective variety, Poincar\'e duality in intersection cohomology implies its invariance under the substitution $q\mapsto q^{-1}$. For example, we have
$$P_{{\bp^n}}(q)=\frac{(-q^\frac{1}{2})^{n+1}-(-q^\frac{1}{2})^{-(n+1)}}{(-q^\frac{1}{2})-(-q^\frac{1}{2})^{-1}}.$$

In particular, we denote by
$$P_{d,m}(q)=P_{M_{d,m}}(q)$$
the shifted Poincar\'e polynomial of the moduli spaces of ordered point configurations.\\[1ex]
Our first main result establishes relations among the Poincar\'e polynomials $P_{d,m}(q)$ under the congruence conditions $m\equiv -1,0,1\bmod d$, respectively:

\begin{thm}\label{r1} For all $d,r\geq 1$, we have a generating function identity in $\mathbb{Q}(q^\frac{1}{2})[[x]]$:

$$\sum_{d\geq 0}\frac{1}{(rd)!}P_{d,{rd+1}}(q)x^d=\exp\left(\sum_{d\geq 1}\frac{1}{(rd)!}P_{\bp^{d-1}}(q)^2P_{{d,rd-1}}(q)x^d\right).$$

Moreover, we have 
$$P_{d,rd}(q)=P_{\bp^{d-1}}(q)\cdot P_{{d,rd-1}}(q).$$
\end{thm}

For coprime $d$ and $m$, we establish a duality
$$M_{d,m}\simeq M_{m-d,m},$$

from which we obtain the following recursion:

\begin{thm}\label{r2} The shifted Poincar\'e polynomials $z_d(q)=P_{{d,2d-1}}(q)$ of the moduli spaces $M_{d,2d-1}$ are given recursively by $z_1(q)=1$ and
$$z_{d+1}(q)=\sum_{k\geq 1}\frac{(2d)!}{k!}\sum_{{a_1+\ldots+a_k=d}\atop{a_1,\ldots,a_k>0}}\prod_{i=1}^k\frac{P_{{\bp}^{a_i-1}}(q)^2\cdot z_{a_i}(q)}{(2a_i)!}.
$$

Equivalently, they are given recursively by $z_1(q)=1$ and
$$z_d(q)=\sum_{{k+l=d}\atop{k,l>0}}z_k(q)z_l(q)[k]^2{{2d-3}\choose{2k-1}}.$$
\end{thm}

Specialization at $q=1$ yields formulas for Gromov-Witten invariants in terms of Euler characteristic $\chi_{(\rm IC)}$ in singular (intersection) cohomology:

\begin{cor}\label{mainresult} For all $d\geq 1$, we have
$$N_{d,0}(\emptyset,(d))=\chi_\mathrm{IC}(M_{d,2d})$$
and
$$N_{d,0}((d),\emptyset)=\chi(M_{d,2d-1})=z_d(1).$$
\end{cor}

Combining the MPS formula \cite{mps} and torus localization \cite{wei}, the Euler characteristic $\chi(M_{d,2d+1})$ is given in a purely combinatorial way in terms of the number of stable spanning trees of a certain full bipartite quiver. More precisely, for a partition $\bp=\sum_{l}lk_l$ of $d$, let $Q(\bp)$ be the full bipartite quiver with $\sum_lk_l$ sinks and $2d+1$ sources. The partition induces a level on the set of sinks by simply choosing $k_l$ sources for each $l$ and assigning them the level $l$. Finally, there is a notion of stability on the set of spanning trees of $Q(\bp)$ where it turns out that a spanning tree is stable if and only if every sink of level $l$ has precisely $2l+1$ neighbours.
\begin{thm}\label{introcoreuler}
We have
\begin{equation*}
\chi(M_{d+1,2d+1})=\chi(M_{d,2d+1}) = \sum_{\bp\vdash d}n(\bp)\prod_{l\geq 1}\frac{(-1)^{k_l(l-1)}l^{k_l(2l-1)}}{k_l!}
\end{equation*}
where $n(\bp)$ is the number of stable spanning trees of $Q(\bp)$.
\end{thm}
It is remarkable that we can use the notion of stable spanning trees in order to obtain an explicit formula for the Euler characteristic which has no Gromov-Witten analogue known to us.
\begin{thm}\label{introeuler}
We have 
$$\chi(M_{d,2d+1})=\frac{(2d)!}{2d+1}\sum_{k=1}^n\frac{(-1)^{d-k}(2d+1)^k}{k!}\sum_{\substack{\sum_{i=1}^k a_i=d\\a_i>0}}\prod_{i=1}^k\frac{a_i^{2a_i-1}}{(2a_i)!}$$
\end{thm}

\begin{rem} Combination of the previous results yields the simple recursion
$$N_{d,0}((d),\emptyset)=\sum_{k+l=d}N_{k,0}((k),\emptyset)N_{l,0}((l),\emptyset)k^2{2d-3\choose 2k-1}$$
for these Gromov-Witten invariants, reminiscent of the Kontsevich recursion for $N_{d,0}$ stated above.\\[1ex]
We expect that, up to shift, the Laurent polynomials $z_d(q)$ are precisely the Block-G\"ottsche invariants \cite{BG} of rational curves in ${\bp}^2$ corresponding to the Gromov-Witten invariants $N_{d,0}((d),\emptyset)$.
\end{rem}

The coefficients of the first few of the polynomials $z_d(q)$ are listed in the following table:

$$\begin{array}{l|l}
d&\mbox{coefficients of }z_d(q)\\ \hline\hline
1&1\\ \hline
2&1\\ \hline
3&1,5,1\\ \hline
4&1,7,29,64,29,7,1\\ \hline
5&1,9,46,175,506,1138,1727,1138,\ldots\\ \hline
6&1,11,67,298,1080,3313,8770,20253,40352,67279,84792,67279\ldots\\ \hline
7&1,13,92,469,1926,6762,20960,58425,148153,344362,735898,\\
&1444761,2591676,4180118,5869613,6735425,5869613,\ldots
\end{array}$$


\section{Quiver moduli methods}\label{quiver}

In this section, we recollect all notions and results on moduli spaces of semistable quiver representations which will be applied in the next section. A general reference is \cite{RModuli}; further references will be given for the respective results.

\subsection{Quiver representations}

Let $Q$ be a finite quiver, consisting of a finite set $Q_0$ of vertices and a finite set of oriented edges (the arrows of $Q$) written $\alpha:i\rightarrow j$. A (complex) representation of $Q$ consists of a tuple $(V_i)_{i\in Q_0}$ of finite dimensional complex vector spaces together with a tuple $(V_\alpha:V_i\rightarrow V_j)_{\alpha:i\rightarrow j}$ of linear maps. A morphism between two such representations $V$ and $W$ consists of a tuple of linear maps $(\varphi_i:V_i\rightarrow W_i)_{i\in Q_0}$ such that all natural squares commute, that is, such that $$\varphi_j\circ V_\alpha=W_\alpha\circ\varphi_i\mbox{ for all }\alpha:i\rightarrow j\mbox{ in }Q.$$
Via componentwise composition, this defines the category ${\bf rep}_\mathbb{C} Q$ of finite dimensional complex representations of $Q$, which is a $\mathbb{C}$-linear abelian category.\\[1ex]
The formal linear combination of dimensions $${\rm\bf dim} V=\sum_{i\in Q_0}(\dim_\mathbb{C}V_i)\cdot i\in\mathbb{N}Q_0$$ is called the dimension vector of $V$. We define the Euler form as the bilinear form on $\mathbb{Z}Q_0$ given by $\langle{\bf d},{\bf e}\rangle=\sum_{i\in Q_0}d_ie_i-\sum_{\alpha:i\rightarrow j}d_ie_j$ for ${\bf d},{\bf e}\in\mathbb{Z}Q_0$. This coincides with the homological Euler form on ${\bf rep}_\mathbb{C}Q$, in the sense that $$\dim_\mathbb{C}{\rm Hom}(V,W)-\dim_\mathbb{C}{\rm Ext}^1(V,W)=\langle{\rm\bf dim} V,{\rm\bf dim W}\rangle,$$ and ${\rm Ext}^{\geq 2}(\_,\_)=0$ on ${\bf rep}_\mathbb{C}Q$.

\subsection{Moduli spaces of quiver representations}

We fix a finite quiver $Q$ without oriented cycles and a dimension vector ${\bf d}\in\mathbb{N}Q_0$ for $Q$. We also fix a form $\Theta=\sum_{i\in Q_0}\Theta_ii^*\in(\mathbb{Z}Q_0)^*$ given by $\Theta({\bf d})=\sum_{i\in Q_0}\Theta_id_i\in\mathbb{Z}$, called a stability for $Q$. This induces a slope function $\mu:\mathbb{N}Q_0\setminus\{0\}\rightarrow\mathbb{Q}$ given by $\mu(V)=\Theta({\rm\bf dim} V)/\dim V$, where $\dim V=\sum_{i\in Q_0}\dim_\mathbb{C}V_i$ denotes the total dimension of $V$.\\[1ex]
We define a representation $V$ of $Q$ to be $\Theta$-semistable (resp.~$\Theta$-stable) if $\mu(U)\leq\mu(V)$ (resp.~$\mu(U)<\mu(V)$) for all non-zero proper subrepresentations $U\subset V$. We call $V$ a $\Theta$-polystable representation if it is isomorphic to a direct sum of stable representations of the same slope (which is then a semisimple object in the abelian subcategory of semistable representations of fixed slope).\\[1ex]
We note that two stabilities $\Theta,\Theta'$ are equivalent, in the sense that the set of stable (resp. semistable, polystable) representations coincide if defined with respect to $\Theta$ and $\Theta'$, respectively, if $$\Theta'=A\Theta+B\dim$$ for suitable $A,B\in\mathbb{Q}$ with $A>0$. In particular, for a given dimension vector ${\bf d}$, every stability $\Theta'$ is equivalent to a stability $\Theta$ such that $\Theta({\bf d})=0$.\\[1ex]
There exists a projective algebraic variety $M_{\bf d}^{\Theta-{\rm sst}}(Q)$ which parametrizes isomorphism classes of $\Theta$-polystable representations of $Q$ of dimension vector ${\bf d}$. It is constructed as follows: we consider the affine space
$$R_{\bf d}(Q)=\bigoplus_{\alpha:i\rightarrow j}{\rm Hom}({\bf C}^{d_i},{\bf C}^{d_j})$$
of all representations of $Q$ of dimension vector ${\bf d}$, on which the group
$$G_{\bf d}=\prod_{i\in Q_0}{\rm GL}_{d_i}({\bf C})$$
acts via base change
$$(g_i)_i\cdot(V_\alpha)_\alpha=(g_j\circ V_\alpha\circ g_i^{-1})_{\alpha:i\rightarrow j},$$
so that $G_{\bf d}$-orbits in $R_{\bf d}(Q)$ naturally correspond to isomorphism classes of representations of $Q$ of dimension vector ${\bf d}$. We denote by $R_{\bf d}^{\Theta-{\rm sst}}(Q)$ (resp.~$R_{\bf d}^{\Theta-{\rm st}}(Q)$) the locus corresponding to $\Theta$-semistable (resp.~$\Theta$-stable) representations. Then
$$M_{\bf d}^{\Theta-{\rm sst}}(Q)=R_{\bf d}^{\Theta-{\rm sst}}(Q)//G_{\bf d}$$
is the quotient parametrizing closed orbits. This moduli space contains a smooth (Zarisiki-)open subset $M_{\bf d}^{\Theta-{\rm st}}(Q)$ consisting of the isomorphism classes of $\Theta$-stable such representations, which can be identified with the geometric quotient of the stable locus by the structure group, that is,
$$M_{\bf d}^{\Theta-{\rm st}}(Q)=R_{\bf d}^{\Theta-{\rm st}}(Q)/G_{\bf d}.$$
Both types of moduli spaces are irreducible of complex dimension $1-\langle{\bf d},{\bf d}\rangle$ if $M_{\bf d}^{\Theta-{\rm st}}(Q)$ is non-empty. We say that ${\bf d}$ is $\Theta$-coprime if $\mu({\bf e})\not=\mu({\bf d})$ for all non-trivial ${\bf e}\leq{\bf d}$. In this case, $\Theta$-stability and $\Theta$-semistability coincide for representations of dimension vector ${\bf d}$, thus $M_{\bf d}^{\Theta-{\rm sst}}(Q)=M_{\bf d}^{\Theta-{\rm st}}(Q)$ is a smooth projective variety. In general, $M_{\bf d}^{\Theta-{\rm sst}}(Q)$ is typically singular.

\subsection{Framed moduli}\label{frame}

Additionally to $Q$, ${\bf d}$ and $\Theta$, we fix ${\bf n}\in\mathbb{N}Q_0$, and let $W_i$ be an $n_i$-dimensional complex vector space, for $i\in Q_0$. We consider ${\bf n}$-framed representations of $Q$, that is, pairs $(V,f)$ consisting of a representation $V$ of $Q$ and a tuple $(f_i:W_i\rightarrow V_i)_{i\in Q_0}$ of linear maps. Such a pair is called $\Theta$-stable if $V$ is $\Theta$-semistable and $\mu(U)<\mu(V)$ whenever $U$ is a non-zero proper subrepresentation of $V$ containing the image of $f$ (that is, ${\rm Im} f_i\subset U_i$ for all $i\in Q_0$). Two ${\bf n}$-framed representations $(V,f)$ and $(V',f')$ are called equivalent if there exists an isomorphism $\varphi:V\rightarrow V'$ intertwining the framings, that is, $f_i'=\varphi_i\circ f_i$ for all $i\in Q_0$.\\[1ex]
There exists a smooth projective variety $M_{{\bf d},{\bf n}}^{\Theta,{\rm fr}}(Q)$ parametrizing equivalence classes of stable ${\bf n}$-framed representations of $Q$ of dimension vector ${\bf d}$. It admits a projective morphism $$\pi:M_{{\bf d},{\bf n}}^{\Theta,{\rm fr}}(Q)\rightarrow M_{\bf d}^{\Theta-{\rm sst}}(Q),$$
given by forgetting the framing datum, whose restriction to the stable locus $M_{\bf d}^{\Theta-{\rm st}}(Q)$ is an \'etale locally trivial fibration with fibre isomorphic to projective space of dimension ${\bf n}\cdot{\bf d}-1$, where ${\bf n}\cdot{\bf d}=\sum_{i\in Q_0}n_id_i$.\\[1ex]
If ${\bf d}$ is $\Theta$-coprime, the map $\pi$ identifies $M_{{\bf d},{\bf n}}^{\Theta,{\rm fr}}(Q)$ with a Zariski-locally trivial $\bp^{{\bf n}\cdot{\bf d}-1}$-bundle over $M_{\bf d}^{\Theta-{\rm sst}}(Q)$ by \cite{RS}.\\[1ex]
For application in the following sections, we recall the construction of the framed moduli space in more detail. Define the framed quiver $\widehat{Q}$ by adding a vertex $i_0$ to $Q_0$, and adding arrows $\beta_{i,k}:i_0\rightarrow i$, for $k=1,\ldots,n_i$ and $i\in Q_0$. Define a dimension vector $\widehat{\bf d}$ for $\widehat{Q}$ by adding an entry $1$ at $i_0$. By choosing bases of all $W_i$, the ${\bf n}$-framed representations $(V,f)$ of $Q$ of dimension vector ${\bf d}$ are naturally identified with representations of $\widehat{Q}$ of dimension vector $\widehat{\bf d}$. We assume that $\Theta$ is normalized such that $\Theta({\bf d})=0$, and extend the stability $\Theta$ for $Q$ to a stability $\widehat{\Theta}$ for $\widehat{Q}$ by adding the entry $1$ at vertex $i_0$. Then we have a natural identification $$M_{{\bf d},{\bf n}}^{\Theta-{\rm fr}}(Q)\simeq M_{\widehat{\bf d}}^{\widehat{\Theta}-{\rm st}}(\widehat{Q}).$$
This follows from \cite[Lemma 3.2]{ER}.

\subsection{Cohomology and Donaldson-Thomas invariants}\label{dt}

For $Q$, ${\bf d}$ and $\Theta$ as above, we define explicit rational functions $P_{\bf d}^\Theta(q)\in\mathbb{Q}(q)$ by
$$P_{\bf d}^\Theta(q)=\sum_{{\bf d}^*}(-1)^{s-1}q^{-\sum_{k\leq l}\langle {\bf d}^l,{\bf d}^k\rangle}\prod_{k=1}^s\prod_{i\in Q_0}\prod_{j=1}^{d_i^k}(1-q^{-j})^{-1},$$
where the sum ranges over all ordered decompositions ${\bf d}={\bf d}^1+\ldots+{\bf d}^s$ of ${\bf d}$ into non-zero dimension vectors such that $$\mu({\bf d}^1+\ldots+{\bf d}^k)>\mu({\bf d})\mbox{ for all }k<s.$$
These functions arise from a resolved Harder-Narasimhan recursion \cite{HNS} and determine the Betti numbers in cohomology of the moduli spaces at least in the $\Theta$-coprime case. Namely, in this case, the normalized Poincar\'e polynomial of $M_{\bf d}^{\Theta-{\rm sst}}(Q)$ is given by
$$P_{M_{\bf d}^{\Theta-{\rm sst}}(Q)}(q)=(q^{-\frac{1}{2}}-q^{\frac{1}{2}})(-q^{\frac{1}{2}})^{\langle{\bf d},{\bf d}\rangle}P_{\bf d}^\Theta(q).$$


In more generality, these rational functions are used to define motivic Donaldson-Thomas invariants as follows:\\[1ex]
For a fixed slope $\mu\in\mathbb{Q}$, define $\Lambda^+_\mu$ as the set of all non-zero dimension vectors of slope $\mu$. Consider the complete commutative local ring $\mathbb{Q}(q^\frac{1}{2})[[\Lambda^+_\mu]]$ with topological basis $x^{\bf d}$ for ${\bf d}\in\Lambda^+_\mu$ and multiplication $x^{\bf d}x^{\bf e}=x^{{\bf d}+{\bf e}}$, with its maximal ideal $\mathfrak{m}$. We consider the plethystic exponential $${\rm Exp}:\mathfrak{m}\rightarrow 1+\mathfrak{m},$$
which is defined as follows:\\[1ex]
For $i\geq 1$, we denote by $\psi_i$ the operator substituting $q$ by $q^i$ and $x^{\bf d}$ by $x^{i{\bf d}}$, and combine them into the operator $\Psi=\sum_{i\geq 1}\frac{1}{i}\psi_i$ on $\mathfrak{m}$. Finally, we define
$${\rm Exp}=\exp\circ\Psi,$$
 which therefore fulfills
$${\rm Exp}(q^ix^{\bf d})=(1-q^ix^{\bf d})^{-1}\mbox{ and }{\rm Exp}(f+g)={\rm Exp}(f)\cdot{\rm Exp}(g).$$

Then we define Donaldson-Thomas invariants ${\rm DT}_{\bf d}^\Theta(q)\in\mathbb{Q}(q^\frac{1}{2})$ by (\cite{KS})

$$1+\sum_{{\bf d}\in\Lambda^+_\mu}(-q^\frac{1}{2})^{\langle{\bf d},{\bf d}\rangle}P_{\bf d}^\Theta(q)x^{\bf d}={\rm Exp}(\frac{1}{q^{-\frac{1}{2}}-q^\frac{1}{2}}\sum_{{\bf d}\in\Lambda^+_\mu}{\rm DT}_{\bf d}^\Theta(q)x^{\bf d}).$$

These invariants ${\rm DT}_d^\Theta(q)$ arise in the following two geometric contexts:\\[1ex]
If the restriction of the Euler from $\langle\_,\_\rangle$ to $\Lambda^+_\mu$ is symmetric and $M_{\bf d}^{\Theta-{\rm sst}}(Q)\not=\emptyset$, then the Donaldson-Thomas invariant determines the intersection Betti numbers of the moduli space $M_{\bf d}^{\Theta-{\rm sst}}(Q)$ by \cite{MeR}:

$$\sum_i \dim {\rm IH}^i(M_{\bf d}^{\Theta-{\rm sst}}(Q),\mathbb{Q})(-q^\frac{1}{2})^i=(-q^\frac{1}{2})^{1-\langle{\bf d},{\bf d}\rangle}{\rm DT}_{\bf d}^\Theta(q).$$

Moreover, the Donaldson-Thomas invariants are related to the cohomology of the framed moduli spaces as follows:

$$1+\sum_{{\bf d}\in\Lambda^+_\mu}P_{M_{{\bf d},{\bf n}}^{\Theta,{\rm fr}}(Q)}(q)(-1)^{{\bf n}\cdot{\bf d}}x^{\bf d}={\rm Exp}(\sum_{\bf d}P_{{\bp}^{{\bf n}\cdot{\bf d}-1}}(q){\rm DT}_{\bf d}^\Theta(q)(-1)^{{\bf n}\cdot{\bf d}}x^{\bf d}).$$

Since this formula is only implicitly contained in \cite{MeR}, we sketch a proof here: we first apply \cite[Theorem 5.2]{ER}, which states that
$$1+\sum_{{\bf d}\in\Lambda^+_\mu}P_{M_{d,n}^{\Theta-{\rm fr}}(Q)}(q)(-q^\frac{1}{2})^{{\bf n}\cdot{\bf d}-\langle{\bf d},{\bf d}\rangle}t^{\bf d}=(1+\sum_{{\bf d}\in\Lambda^+_\mu}P_{\bf d}^\Theta(q)t^{\bf d})^{-1}\cdot(1+\sum_{{\bf d}\in\Lambda^+_\mu}q^{{\bf n}\cdot{\bf d}}P_{\bf d}^\Theta(q)t^{\bf d}),$$
where the formal variables $t^{\bf d}$ multiply by
$$t^{\bf d}\cdot t^{\bf e}=q^{-\langle{\bf d},{\bf e}\rangle}t^{{\bf d}+{\bf e}},$$
and where we adjusted to the present definition of the (shifted) Poincar\'e polynomials, noting vanishing of odd cohomology of $M_{{\bf d},{\bf n}}^{\Theta-{\rm fr}}(Q)$. Since the restriction of the Euler form to $\Lambda^+_\mu$ is assumed to be symmetric, we can substitute $t^{\bf d}$ by $(q^\frac{1}{2})^{\langle{\bf d},{\bf d}\rangle}x^{\bf d}$ to derive the following equality in $\mathbb{Q}(q^\frac{1}{2})[[\Lambda^+_\mu]]$:
$$1+\sum_{{\bf d}\in\Lambda^+_\mu}P_{M_{d,n}^{\Theta-{\rm fr}}(Q)}(q)(-q^\frac{1}{2})^{{\bf n}\cdot{\bf d}}x^{\bf d}=\frac{1+\sum_{{\bf d}\in\Lambda^+_\mu}q^{{\bf n}\cdot{\bf d}}P_{\bf d}^\Theta(q)(-q^\frac{1}{2})^{\langle{\bf d},{\bf d}\rangle}x^{\bf d})}{1+\sum_{{\bf d}\in\Lambda^+_\mu}P_{\bf d}^\Theta(q)(-q^\frac{1}{2})^{\langle{\bf d},{\bf d}\rangle}x^{\bf d})}.$$
Using the defining equation for Donaldson-Thomas invariants, this simplifies to
$$1+\sum_{{\bf d}\in\Lambda^+_\mu}P_{M_{d,n}^{\Theta-{\rm fr}}(Q)}(q)(-q^\frac{1}{2})^{{\bf n}\cdot{\bf d}}x^{\bf d}={\rm Exp}(\sum_{{\bf d}\in\Lambda^+_\mu}\frac{q^{{\bf n}\cdot{\bf d}}-1}{q^{-\frac{1}{2}}-q^\frac{1}{2}}{\rm DT}_{\bf d}^\Theta(q)x^{\bf d}).$$
We identify
$$\frac{q^N-1}{q^{-\frac{1}{2}}-q^\frac{1}{2}}=(-q^\frac{1}{2})^NP_{\bp^{N-1}}(q)$$
and note that formation of ${\rm Exp}$ is compatible with the replacement of $q^{\frac{1}{2}{\bf n}\cdot{\bf d}}x^{\bf d}$ by $x^{\bf d}$. Then the previous equality reduces to the claimed one.

\subsection{Small resolutions}\label{small}

Assume that ${\bf d}$ is an indivisible dimension vector, that $\Theta$ is normalized such that $\Theta({\bf d})=0$, that $M_{\bf d}^{\Theta-{\rm st}}(Q)\not=\emptyset$, and that the restriction of $\langle\_,\_\rangle$ to $\Lambda^+_0$ is symmetric.\\[1ex]
Since ${\bf d}$ is indivisible, there exists a stability $\eta$ such that $\eta({\bf d})=0$ and $\eta({\bf e})\not=0$ whenever $0\not={\bf e}\leq{\bf d}$ is such that $\Theta({\bf e})=0$. We choose $C\in\mathbb{N}$ such that
$$C>\max(\max(\eta({\bf e}\, :\, {\bf e}\leq{\bf d},\, \Theta(n{\bf e})<0),\max(-\eta({\bf e})\, :\, {\bf e}\leq{\bf d},\, \Theta({\bf e})>0)))$$
and define $\Theta'=C\Theta+\eta$. Then by \cite[Theorem 4.3, Theorem 5.1]{RSm}, there exists a small resolution
$$M_{\bf d}^{\Theta'-{\rm sst}}(Q)\rightarrow M_{\bf d}^{\Theta-{\rm sst}}(Q)$$
and
$$P_{M_{\bf d}^{\Theta'-{\rm sst}}(Q)}(q)={\rm DT}_{\bf d}^\Theta(q).$$

\subsection{MPS degeneration formula for moduli spaces}\label{mpsform}
Fix a quiver $Q$ and a linear $\Theta$ as above. For a vertex $r\in Q_0$, we denote by $A_r\subseteq Q_1$ the set of arrows $\alpha$ such that $r=h(\alpha)$ or $r=t(\alpha)$, i.e. where $r$ is the head or tail of $\alpha$. Moreover, let $Q(r)$ be the quiver which has vertices
\[Q(r)_0=Q_0\backslash \{r\}\cup\{r_{l,n}\mid (l,n)\in\Nn_+^2\}\]
and arrows
\begin{eqnarray*}
Q(r)_1&=&Q_1\backslash A_r\cup\{\alpha_1,\ldots,\alpha_l:p\rightarrow r_{l,n}\mid \alpha:p\rightarrow r,\,l,n\in\Nn_+\}\\&&\cup\{\alpha_1,\ldots,\alpha_l:r_{l,n}\rightarrow q\mid\alpha:r\rightarrow q,\,n\in\Nn_+\}.
\end{eqnarray*}
We define a level $l:Q(r)_0\rightarrow\mathbb{N}$ by $$l(q)=\begin{cases}l\text{ if }q=r_{l,n}\\1\text{ if }q\in Q_0\backslash\{r\}\end{cases}.$$ 

This defines a linear form $\Theta^r\in\Hom(\Zn Q(r)_0,\Zn)$ by $(\Theta^r)_q=\Theta_q$ for all $q\not=r$  and $(\Theta^r)_{r_{l,n}}=l\Theta_r$ for all $l,n\geq 1$. We consider the slope $\mu=\Theta^r/\kappa$ where $\kappa$ is defined by $\kappa(\bd)=\sum_{q\in Q(r)_0}l(q)\bd_q$ for $\bd\in Q(r)_0$. We denote the corresponding moduli space of stable representations by $M_\bd^{\Theta^r-\st}(Q(r))$. 

If we fix a dimension vector $\bd\in\Nn Q_0$ and a (weighted) partition $\bd_r=\sum_{l=1}^tlk_l$, this induces a dimension vector $\bar{\bd}$ of $Q(r)$ in the following way: we set $\bar{\bd}_q=\bd_q$ for all $q\neq r$ and $\bar \bd{r_{l,n}}=1$ for $1\leq l\leq t$ and $1\leq n\leq k_l$ and $\bar \bd_{r_{l,n}}=0$ otherwise. If we think of a dimension vector of $Q(r)$, it is convenient to think of a tuple $\bd(k_*):=(\bd,k_*)$ where $k_*\vdash \bd_r$ is a weighted partition. We call $Q(r)$ the MPS-quiver of $Q$ with respect to $r$. Clearly, we can inductively apply this construction to all vertices of $Q$ and obtain the full MPS-quiver of $Q$. 

Now we can formulate the following result concerning the Euler characteristic of moduli spaces of stable representations, see \cite[Appendix D]{mps} and also \cite[Sections 3.2, 3.3]{rsw} for a more general setting:
\begin{thm}\label{mpsthm} If $\bd$ is $\Theta$-coprime, we have

and
$$\chi(M^{\Theta-\st}_\bd(Q))=\sum_{k_*\vdash \bd_r}\prod_{l\geq 1}\frac{1}{k_l!}\left(\frac{(-1)^{l-1}}{l^2}\right)^{k_l}\chi(M_{\bd(k_*)}^{\Theta^r-\st}(Q(r)))$$
where the sum ranges over all weighted partitions of $\bd_r$.
\end{thm}

\subsection{Localization for thin dimension vectors}\label{loc}
Fix $Q$ and $\Theta$ as above. In order to derive a purely combinatorial description of $\chi(M^{\Theta-\st}_\bd(Q))$, we combine the second part of Theorem \ref{mpsthm} with the localization theorem which is particularly powerful in the case of thin dimension vectors, i.e. if $\bd_q\in\{0,1\}$ for every $q\in Q_0$. We assume that $\bd$ is thin and supported at $Q_0$. In this case every spanning tree of $Q$, i.e. every subquiver $T$ of $Q$ with $T_0=Q_0$ such that the underlying graph is a connected tree, naturally defines a representation of dimension $\bd$ by assigning the identity map to each arrow of the tree and the zero map to the remaining arrows of $Q$. We use this observation to introduce the notion of stable spanning trees. In this case the localization theorem \cite[Corollary 3.14]{wei} can be formulated as follows.
\begin{thm}\label{locforthin}
If $\bd$ is a thin dimension vector, the Euler characteristic $\chi(M^{\Theta-\mathrm{st}}_\bd(Q))$ is given by the number of stable spanning trees of $Q$. 
\end{thm}

\section{Techniques for moduli spaces of point configurations}\label{point}

\subsection{Quiver moduli setup for moduli spaces of point configurations}

As before, we fix nonnegative integers $d$ and $m$. We consider the quiver $Q_m$ with vertices $i_1,\ldots,i_m$ and $j$ and arrows $i_k\rightarrow j$ for all $k=1,\ldots,m$, called the $m$-subspace quiver. For a subset $I\subset[m]$ and a nonnegative integer $e$, we denote by ${\bf d}(I,e)$ the dimension vector $${\bf d}(I,e)=\sum_{k\in I}i_k+e\cdot j\in\mathbb{N}(Q_m)_0$$ for $Q_m$. In particular, we denote ${\bf d}={\bf d}([m],d)$.\\[1ex]
The variety $R_{\bf d}(Q_m)$ is then isomorphic to the affine space $(\mathbb{C}^d)^m$, and the group $G_{\bf d}=(\mathbb{C}^*)^m\times{\rm GL}_d(\mathbb{C})$ acts on it via
$$((t_1,\ldots,t_m),g)\dot(v_1,\ldots,v_m)=(t_1^{-1}gv_1,\ldots,t_m^{-1}gv_m)$$
for $t_1,\ldots,t_m\in\mathbb{C}^*$, $g\in{\rm GL}_d(\mathbb{C})$ and $v_1,\ldots,v_m\in\mathbb{C}^d$.\\[1ex]
Writing $d=g\overline{d}$ and $m=g\overline{m}$ for $g=\gcd(d,m)$, we define a stability $\Theta$ for $Q_m$ by $$\Theta=\overline{d}\sum_{k=1}^mi_k^*-\overline{m}j^*,$$
thus $\Theta({\bf d})=0$.\\[1ex]
Translating the definition of moduli spaces of semistable quiver representations to the present setup of the quiver $Q_m$, the dimension vector ${\bf d}$ and the stability $\Theta$, we see that $$M_{d,m}\simeq M_{\bf d}^{\Theta-{\rm sst}}(Q_m)$$
is the moduli space of point configurations considered above.
We consider the dimension vector ${\bf n}=j$ and also make use of the framed moduli space
$$M_{d,m}^{\rm fr}=M_{{\bf d},{\bf n}}^{\Theta,{\rm fr}}(Q).$$
It parametrizes semistable point configurations $(p_1,\ldots,p_m)$, together with a point $p_0\in\mathbb{P}^{d-1}$ such that $p_0\not\in U_I$ whenever $\dim U_I= \frac{d}{m}|I|-1$ for a non-empty proper subset $I$, up to projective equivalence.




\subsection{Duality}

In this section, we prove the following duality:

\begin{thm} For coprime $d$ and $m$ such that $d<m$, we have
$$M_{d,m}\simeq M_{m-d,m}.$$
\end{thm}



Denote by ${\bf d'}$ the dimension vector $\sum_{k=1}^mi_k+(m-d)j$ for $Q_m$. For nonnegative integers $a$ and $b$, denote by $M_{a\times b}(\mathbb{C})$ the affine space of complex $a\times b$-matrices, and by $M_{a,b}(\mathbb{C})'$ the open subset of matrices of maximal rank. We identify $R_{\bf d}(Q_m)$ with $M_{d\times m}(\mathbb{C})$ by associating to a matrix $A$ its tuple of columns; in the same way, we identify $R_{\bf d'}(Q_m)$ with $M_{(m-d)\times m}(\mathbb{C})$.\\[1ex]
Inside $M_{d\times m}(\mathbb{C})'\times M_{(m-d)\times m}(\mathbb{C})'$, we consider the closed subvariety $Z$ of pairs $(A,B)$ of matrices of maximal rank such that $AB^T=0$. The group ${\rm GL}_d(\mathbb{C})\times(\mathbb{C}^*)^m\times{\rm GL}_{m-d}(\mathbb{C})$ naturally acts on $Z$, and we have natural projection maps $$M_{d\times m}(\mathbb{C})'\stackrel{p_1}{\leftarrow} Z\stackrel{p_2}{\rightarrow}M_{(m-d)\times m}(\mathbb{C})'.$$
Since the equation $AB^T=0$ means that the columns of $B^T$ form $m-d$ linearly independent solutions to the homogeneous linear equation $Ax=0$, every fibre of $p_1$ is a single ${\rm GL}_{m-d}(\mathbb{C})$-orbit. Moreover, on the open subset of $M_{d\times m}(\mathbb{C})'$ where a fixed $d\times d$-submatrix is invertible, we can use this submatrix to construct, polynomially in the entries of $A$, a matrix $B$ such that $(A,B)\in Z$. These two facts together yield identifications of geometric quotients
$$M_{d\times m}(\mathbb{C})'\simeq Z/{\rm GL}_{m-d}(\mathbb{C})\mbox{ and similarly }M_{(m-d)\times m}(\mathbb{C})'\simeq Z/{\rm GL}_d(\mathbb{C}).$$
Using the identification of $M_{d\times m}(\mathbb{C})$ with $R_{\bf d}(Q_m)$, we call a matrix $A$ semistable if its tuple of columns defines a $\Theta$-semistable point in $R_{\bf d}(Q_m)$. Using the definition of semistability, this translates to the following property of $A$:\\[1ex]
For every subset $I\subset[m]$, let $A_I$ be the submatrix consisting of the columns $A_i$ for $i\in I$. Then $A$ is semistable if and only if for all such $I$ we have ${\rm rk}(A_I)\geq \frac{d}{m}|I|$. Similarly, a matrix $B\in M_{(m-d)\times m}(\mathbb{C})$ is semistable via the identification of $M_{(m-d)\times m}(\mathbb{C})$ with $R_{\bf d'}(Q_m)$ if and only if ${\rm rk}(B_I)\geq\frac{m-d}{m}|I|$ for all such $I$.\\[1ex]
We claim that, for $(A,B)\in Z$, semistability of $A$ and of $B$ are equivalent. We need the following easy linear algebra lemma whose proof is left to the reader:

\begin{lem} Given a square of vector spaces and linear maps
$$\begin{array}{lcr}V_1&\stackrel{f}{\rightarrow}&V_2\\ g\downarrow&&\downarrow h\\ V_3&\stackrel{i}{\rightarrow}&V_4\end{array}$$
inducing a short exact sequence
$$0\rightarrow V_1\stackrel{[{f\atop g}]}{\rightarrow}V_2\oplus V_3\stackrel{[h\, i]}{\rightarrow}V_4\rightarrow 0,$$
the map $g$ induces an isomorphism $g:{\rm Ker}(f)\stackrel{\simeq}{\rightarrow}{\rm Ker}(i)$.
\end{lem}

Now suppose that $A$ is semistable and that $(A,B)\in Z$, thus $B$ has maximal rank and $AB^T=0$. Let $I$ be a subset of $[m]$; we have to prove that ${\rm rk}(B_I)\geq\frac{m-d}{m}|I|$. We denote by $J$ the complement of $I$ and consider the square of linear maps
$$\begin{array}{lcr}\mathbb{C}^{m-d}&\stackrel{B_I^T}{\rightarrow}&\mathbb{C}^{|I|}\\ B_J^T\downarrow&&\downarrow A_I\\ \mathbb{C}^{|J|}&\stackrel{A_J}{\rightarrow}&\mathbb{C}^d.\end{array}$$
By assumption the condition of the previous lemma is satisfied, so that $${\rm Ker}(B_I^T)\simeq{\rm Ker}(A_J).$$
This allows us to estimate:
$${\rm rk}(B_I)={\rm rk}(B_I^T)=m-d-\dim{\rm Ker}(B_I^T)=$$
$$=m-d-\dim{\rm Ker}(A_J)=m-d-|J|+{\rm rk}(A_J)\geq$$
$$\geq m-d-|J|+\frac{d}{m}|J|=m-d-\frac{m-d}{m}|J|=\frac{m-d}{m}(m-|J|)=$$
$$=\frac{m-d}{m}|I|,$$
as desired. Dually, semistabilty of $B$ implies semistability of $A$.\\[1ex]
We denote by $Z^{\rm sst}$ the open subset of pairs $(A,B)\in Z$ where $A$, or equivalently $B$, is semistable. Since by coprimality of $d$ and $m$, the moduli spaces $M_{d,m}$ and $M_{m-d,m}$ are geometric quotients of the corresponding (semi-)stable loci, we derive the claimed duality by the following series of identifications of geometric quotients:
$$M_{d,m}=R_{\bf d}^{\Theta_0-{\rm sst}}(Q_m)/({\rm GL}_d(\mathbb{C})\times(\mathbb{C}^*)^m)\simeq$$
$$\simeq Z^{\rm sst}/({\rm GL}_{d}(\mathbb{C})\times(\mathbb{C}^*)^m\times{\rm GL}_{m-d}(\mathbb{C}))$$
$$\simeq R_{\bf d'}^{\Theta-{\rm sst}}(Q_m)/((\mathbb{C}^*)^m\times{\rm GL}_{m-d}(\mathbb{C}))=M_{m-d,m},$$
proving the claimed duality.

\subsection{Derivation of a recursion}

In this section, we establish a recursion relating the shifted Poincar\'e polynomials of the framed moduli spaces of point configurations $M_{d,m}^{\rm fr}$ with those of the ordinary moduli spaces $M_{d,m}$.\\[1ex]
A direct computation easily shows that the Euler form of $Q_m$ is symmetric on $\Lambda^+_0$, which by Section \ref{dt} yields
$${\rm DT}_{\bf d}^\Theta(Q_m)=P_{d,m}(q).$$
Next we determine those dimension vectors ${\bf d}(I,e)$ which belong to the lattice $\Lambda^+_0$. Again we write $d=g\overline{d}$ and $m=g\overline{m}$ for $g=\gcd(d,m)$. So we have $\gcd(\overline{d},\overline{m})=1$, and thus $$0=\Theta({\bf d}(I,e)))=\overline{m}e-\overline{d}|I|\mbox{ if and only if }|I|=a\overline{m}\mbox{ and }e=a\overline{d}\mbox{ for }0\leq a\leq g.$$
For such a pair $(I,e)$, we have $$M_{{\bf d}(I,e)}^{\Theta-{\rm sst}}(Q_m)\simeq M_{a\overline{m},a\overline{d}}$$ and thus $${\rm DT}_{{\bf d}(I,e)}^{\Theta}(q)=P_{a\overline{d},a\overline{m}}(q).$$

Our aim is to prove the following recursion:

\begin{thm}\label{mainrecursion} In the notation above, we have
$$P_{M_{d,m}^{\rm fr}}(q)=\sum_{k\geq 1}\frac{m!}{k!}\sum_{a_1+\ldots+a_k=g}\prod_{i=1}^k\frac{P_{{\bp}^{a_i\overline{d}-1}}(q)P_{a_i\overline{d},a_i\overline{m}}(q)}{(a_i\overline{m})!}.$$
\end{thm}

Again we fix the framing datum ${\bf n}=j$ for $Q_m$. We would like to compare the $x^{\bf d}$-terms on both sides of the equation
$$1+\sum_{{\bf d}\in\Lambda^+_0}P_{M_{{\bf d},{\bf n}}^{\Theta,{\rm fr}}(Q)}(q)(-1)^{{\bf n}\cdot{\bf d}}x^{\bf d}={\rm Exp}(\sum_{\bf d}P_{{\bp}^{{\bf n}\cdot{\bf d}-1}}(q){\rm DT}_{\bf d}^\Theta(q)(-1)^{{\bf n}\cdot{\bf d}}x^{\bf d})$$
from Section \ref{dt}. On the left hand side, this term is obviously $(-1)^{d}P_{M_{d,m}^{\rm fr}}(q)$.\\[1ex]
On the right hand side, we can, and will, work modulo $x_{i_k}^2$ for $k=1,\ldots,m$. For $I\subset[m]$, we denote by $x_I$ the monomial $\prod_{k\in I}x_{i_k}$. Then, by the previous consideration of dimension vectors in $\Lambda^+_0$, the right hand side simplifies to
$${\rm Exp}(\sum_{a=1}^g\sum_{{I\subset[m]}\atop{|I|=a\overline{m}}}P_{{\bp}^{a\overline{d}-1}}(q){\rm DT}_{{\bf d}(I,a\overline{d})}^{\Theta}(q)(-1)^{a\overline{d}}x_Ix_j^{a\overline{d}}).$$
By definition, ${\rm Exp}$ is calculated by first applying $\Psi$, then $\exp$. We note that $\Psi$ acts as the identity modulo $x_{i_k}^2$. Thus it suffices to apply $\exp$ modulo $x_{i_k}^2$, easily yielding
$$\sum_{k\geq 1}\frac{1}{k!}\sum_{a_1+\ldots+a_k\leq g}\sum_{I_1,\ldots,I_k}\prod_{i=1}^k(P_{{\bp}^{a_i\overline{d}-1}}(q){\rm DT}_{{\bf d}(I_i,a_i\overline{d})}^{\Theta_0}(q)(-1)^{a_i\overline{d}}x_{I_i}x_j^{a_i\overline{d}}),$$
where the inner sum ranges over tuples $(I_1,\ldots,I_k)$ of pairwise disjoint subset of $[m]$ with cardinality $|I_i|=a_i\overline{m}$ for all $i$. Thus the coefficient of $x^{\bf d}=x_{[m]}x_j^d$ equals
$$\sum_{k\geq 1}\frac{1}{k!}\sum_{a_1+\ldots+a_k=g}\sum_{I_1\cup\ldots\cup I_k=[m]}\prod_{i=1}^k(P_{{\bp}^{a_i\overline{d}-1}}(q){\rm DT}_{{\bf d}(I_i,a_i\overline{d})}^{\Theta_0}(q)(-1)^{a_i\overline{d}}),$$
the sum ranging over all partitions $[m]=I_1\cup\ldots\cup I_k$ into pairwise disjoint subsets of cardinality $|I_i|=a_i\overline{m}$ for all $i$. For a fixed ordered decomposition $g=a_1+\ldots+a_k$, the number of such partitions equals
$$\frac{m!}{(a_1\overline{m})!\cdot\ldots\cdot(a_k\overline{m})!}.$$ so the previous term simplifies to
$$(-1)^{d}\sum_{k\geq 1}\frac{1}{k!}\sum_{a_1+\ldots+a_k=g}\frac{m!}{(a_1\overline{m})!\cdot\ldots\cdot(a_k\overline{m})!}\prod_{i=1}^k(P_{{\bp}^{a_i\overline{d}-1}}(q)P_{a_i\overline{d},a_i\overline{m}}(q)).$$
The claimed identity follows.

\subsection{From framed to unframed moduli spaces}

In this section we derive the following three relations between framed and unframed moduli spaces for special $d$ and $m$.
\begin{pro}\label{propfrufr} We have the following relations between framed and unframed moduli spaces:
\begin{enumerate}
\item If $d$ and $m$ are coprime, the natural map $M_{d,m}^{\rm fr}\rightarrow M_{d,m}$ turns $M_{d,m}^{\rm fr}$ into a Zariski-locally trivial ${\bp}^{d-1}$-bundle over $M_{d,m}$.
\item If $d$ divides $m$, the framed moduli space $M_{d,m}^{\rm fr}$ is isomorphic to $M_{d,m+1}$.
\item If $d$ divides $m$, there exists a small resolution $M_{d,m-1}^{\rm fr}\rightarrow M_{d,m}$.
\end{enumerate}
\end{pro}

The first statement is a direct consequence of the general result on framed moduli spaces.\\[1ex]
For the second statement, we first note that in case $m=rd$, the moduli space $M_{d,m}$ is defined using the stability $$\Theta=\sum_{k=1}^{rd}i_{k}^*-rj^*.$$
The framing construction from Section \ref{frame} then yields an isomorphism
$$M_{d,m}^{\rm fr}\simeq M_{\widehat{\bf d}}^{\widehat{\Theta}-{\rm sst}}(Q_{m+1})$$ for the stability $$\widehat{\Theta}=\sum_{k=1}^{m+1}i_k^*-rj^*$$ for $Q_{m+1}$. We replace $\widehat{\Theta}$ by the equivalent stability $$\Theta'=(d+\frac{1}{r+1})\widehat{\Theta}-\frac{1}{r+1}\dim$$ and arrive at the stability $$\Theta'=\sum_{k=1}^{m+1}i_k^*-(m+1)j^*$$ for $Q_{m+1}$, which is used to define the moduli space $M_{d,m+1}$, as desired.\\[1ex]
For the third statement, we start with the stability $$\Theta=d\sum_{k=1}^{m-1}i_k^*-(m-1)j^*$$ for $Q_{m-1}$ which is used to define the moduli space $M_{d,m-1}$. The framing construction yields an isomorphism $$M_{d,m-1}^{\rm fr}\simeq M_{\widehat{\bf d}}^{\widehat{\Theta}-{\rm sst}}(Q_m)$$ for the stability
$$\widehat{\Theta}=i_1+d\sum_{k=2}^m i_k^*-(m-1)j^*$$
for $Q_m$. Writing again $m=rd$, we replace $\widehat{\Theta}$ by the equivalent stability $$\Theta'=(r+1)d\widehat{\Theta}-\dim.$$
For the stabilities $$\Theta''=\sum_{k=1}^mi_k^*-rj^*\mbox{ and }\eta=-(r+1)d(d-1)i_1^*+(r+1)(d-1)j^*$$ for $Q_m$ and $C=(r+1)d^2-1$, a direct computation shows that $$\Theta'=C\Theta''+\eta.$$
Evaluating $\eta$ at all dimension vectors ${\bf d}(I,e)$, we easily derive the estimate $$|\eta({\bf d}(I,e)|\leq (r+1)d(d-1)<C$$ for all $I$ and $e$. From Section \ref{small} it follows that $M_{\bf d}^{\Theta'-{\rm sst}}(Q_m)$ is a small resolution of $M_{\bf d}^{\Theta''-{\rm sst}}(Q)$, which in turn is isomorphic to $M_{d,m}$. But on the other hand $M_{\bf d}^{\Theta'-{\rm sst}}(Q_m)$ is isomorphic to $M_{d,m-1}^{\rm fr}$, as claimed.

\section{Applications and proofs of the main results}\label{app}

In this section we prove most of the statements of Section \ref{results}, as well as giving two further applications of the above methods.

\subsection{Recursion for unframed moduli spaces}

We start with the main recursion Theorem \ref{mainrecursion} applied to the case $m=rd$, giving

$$P_{M_{d,rd}^{\rm fr}}(q)=\sum_{k\geq 1}\frac{(rd)!}{k!}\sum_{a_1+\ldots+a_k=d}\prod_{i=1}^k\frac{P_{{\bp}^{a_i\overline{d}-1}}(q){\rm DT}_{a_i,ra_i}(q)}{(ra_i)!}.$$

We use Proposition \ref{propfrufr} to rewrite
$$P_{M_{d,rd}^{\rm fr}}(q)=P_{{d,rd+1}}(q)$$
and
$${P}_{d,rd}(q)=P_{M_{d,rd-1}^{\rm fr}}(q)=P_{\bp^{d-1}}(q)\cdot P_{{d,rd-1}}(q).$$
Then the above recursion reads
$$P_{{d,rd+1}}(q)=\sum_{k\geq 1}\sum_{a_1+\ldots+a_k=d}\prod_{i=1}^k\frac{P_{\bp^{a_i-1}}(q)^2P_{{a_i,ra_i-1}}(q)}{(ra_i)!}.$$
We need a general lemma on recursions of this type:

\begin{lem} For sequences $(a_n)_{n\geq 1}$ and $(b_n)_{n\geq 1}$ in a field $K$ of characteristic $0$ the following are equivalent:
\begin{enumerate}
\item The sequences are related by the recursion
$$a_n=\sum_{s=1}^n\frac{1}{s!}\sum_{n=n_1+\ldots+n_s}b_{n_1}\cdot\ldots\cdot b_{n_s}$$
for all $n\geq 1$,
\item The generating series $F(t)=1+\sum_{n\geq 1}a_nt^n$ and $G(t)=\sum_{n\geq 1}b_nt^n$ in $K[[t]]$ are related by
$$F(t)=\exp(G(t)),$$
\item The sequences are related by the recursion
$$na_n=\sum_{k+l=n}kb_ka_l$$
for all $n\geq 1$.
\end{enumerate}
\end{lem}

The first relation in the statement of the lemma is equivalent to the second by expanding the exponential series. Taking the differential $\frac{d}{dt}$ on both sides yields
$$\frac{d}{dt}F(t)=\frac{d}{dt}G(t)\cdot F(t),$$
which, after comparing coefficients, is equivalent to the third recursion.\\[1ex]
We apply this lemma to the sequences
$$a_d=\frac{1}{(rd)^!}P_{{d,rd+1}}(q),\;\;\; b_d=\frac{1}{(rd)!}P_{\bp^{d-1}}(q)^2P_{{d,rd-1}}(q)$$
and obtain the recursion
$$dP_{{d,rd+1}}(q)=\sum_{e+f=d}eP_{\bp^{e-1}}(q)^2P_{{e,re-1}}(q)P_{{f,rf+1}}(q){rd\choose re}$$
and the generating function equation
$$\sum_{d\geq 0}\frac{1}{(rd)!}P_{d_{rd+1}}(q)x^d=\exp(\sum_{d\geq 1}\frac{1}{(rd)!}P_{\bp^{d-1}}(q)^2P_{{d,rd-1}}(q)x^d),$$
proving Theorem \ref{r1}.

\subsection{Moduli spaces of point configurations in the projective line}

As a side remark, we will now derive a simple recursion for the shifted Poincar\'e polynomials of the moduli spaces $M_{2,m}$ of point configurations in the projective line, which can be compared to the formulas of \cite{MFK}.\\[1ex]
We consider the special case $d=2$ and work out the recursion of Theorem \ref{mainrecursion}:

$$P_{{2,2r+1}}(q)=(2r)!\frac{P_{\bp^1}(q)^2P_{2,2r-1}(q)}{(2r)!}+\frac{(2r)!}{2}\frac{P_{{1,r-1}}(q)^2}{(r!)^2}.$$
This easily simplifies to the recursion

$$P_{{2,2r+1}}(q)=P_{\bp^1}(q)^2P_{{2,2r+1}}(q)+\frac{1}{2}{2r\choose r}$$
for moduli spaces of an odd number of points. The even case reduces to this via
$$P_{2,2r}(q)=P_{\bp^1}(q)\cdot P_{2,2r-1}(q).$$

\subsection{Application to Gromov-Witten invariants}

To derive Theorem \ref{r2}, we now consider the special case $r=2$ of Theorem \ref{mainrecursion}. We then use duality to identify

$$M_{d,2d+1}\simeq M_{d+1,2d+1}=M_{d+1,2(d+1)-1}.$$

Thus we find

$$P_{{d+1,2(d+1)-1}}(q)=\sum_{k\geq 1}\sum_{a_1+\ldots+a_k=d}\prod_{i=1}^k\frac{P_{\bp^{a_i-1}}(q)^2P_{{a_i,2a_i-1}}(q)}{(2a_i)!},$$
and Theorem \ref{r2} is proved. Comparing the specialization of this formula at $q=1$ to the formulas \cite[Corollary 4.8, Theorem 4.9]{fm} for the Gromov-Witten invariants $N_{d,0}((d),\emptyset)$ and $N_{d,0}(\emptyset,(d))$ (which will also be considered in Section \ref{comb}), we arrive at Corollary \ref{mainresult}.

\subsection{An identity between GW invariants}

Following \cite{gps}, we consider the following Gromov-Witten invariants:\\[1ex]
We fix coprime $a,b\in\Nn^+$ and consider the weighted projective plane (resp. toric surface) $X_{a,b}=(\Cn^3\backslash\{0\})/\Cn^\ast$ for the action defined by $t\cdot(x,y,z):=(t^ax,t^by,tz)$. We denote by $X_{a,b}^o$ the open subset obtained when removing the three toric fixed points and by $D_1^o,D_2^o$ and $D_{\mathrm{out}}^o$ the toric divisors without these fixed points.

With every pair of ordered partitions $\bp_1,\bp_2$ of lengths $l_1,l_2$ we can associate a Gromov-Witten invariant $N_{a,b}[({\bf P}_1,{\bf P}_2)]$, see \cite[Section 0.4]{gps}.
Heuristically, this invariant may be viewed as the ``number''  of rational curves in $X_{a,b}$ intersecting $l_i$ distinct fixed points of $D_i^o$ with multiplicities given by the $p_{i,l}$ for $i=1,2$, and intersecting $D_{\rm out}^o$ in one unspecific point with multiplicity one.

By the Refined Gromov-Witten/Kronecker correspondence \cite[Corollary 9.1]{rw} together with Theorem \ref{mainresult}, we get an identity of Gromov-Witten invariants:
\begin{cor}
We have $N_{d,0}((d),\emptyset)=N_{(d,2d-1)}[(d,1^{2d-1})]$.
\end{cor}

\section{Combinatorics}\label{comb}

Using the machinery of tropical geometry on the Gromov-Witten side and the MPS formula together with the localization theorem for quiver moduli on the point configuration side, one obtains an enumerative problem on both sides. We introduce and discuss the combinatorial objects which are to be enumerated, namely labeled floor diagrams and stable trees, and discuss similarities and differences.

\subsection{Labeled floor diagrams}
\noindent We consider (markings of) labeled floor diagrams. We modify the definition in comparison to \cite{fm} using the language of quivers as it turns out to be useful for our purposes. We fix natural number $d>0$ and $g$. . 
\begin{defi}A (connected) labeled floor diagram $\mathcal D$ of degree $d$ and genus $g$ is a quiver with $d$ vertices $Q_0=\{1,\ldots,d\}$ and $d+g-1$ arrows $Q_1$, which are weighted by a weight function $w:Q_1\to\Nn^+$, such that
\begin{enumerate}
\item for each arrow $\alpha:p\to q$, we have $p<q$;
\item for each vertex $q\in Q_0$, we have
\[\mathrm{div}(q):=\sum_{\alpha\in Q_1,t(\alpha)=q}w(\alpha)-\sum_{\alpha\in Q_1,h(\alpha)=q}w(\alpha)\leq 1.\]
\end{enumerate}
The multiplicity of $\mathcal D$ is defined by
\[\mu(\mathcal D)=\prod_{\alpha\in Q_1}(w(\alpha))^2.\]
\end{defi}
Note that a labeled floor diagram has no oriented cycles and loops. For a fixed labeled floor diagram $\mathcal D$ of degree $d$ and genus $g$, we introduce the notion of markings. Again we modify the definition of markings defined in \cite[Section 1]{fm} for our purposes. We fix two partitions (possibly empty) $\lambda$ and $\rho$ such that $|\lambda|+|\rho|=d$. A $(\lambda,\rho)$-marking of a labeled floor diagram $\mathcal D$ is obtained as follows. 

In a first step, we add vertices $V:=\{v_1,\ldots, v_{l(\lambda)}\}$ and $W:=\{w_1,\ldots,w_{l(\rho)}\}$. Moreover, we add weighted arrows from the original vertices to the new vertices in such a way that, from every original vertex $v$, there emanate altogether $(1-\mathrm{div}(v))$ arrows (taking weights into account). Moreover, there is exactly one (weighted) arrow pointing to each new vertex. Finally, we want the weights of the arrows pointing to the vertices $v_1,\ldots, v_{l(\lambda)}$ to coincide with the partition $\lambda$ and those pointing to the vertices $w_1,\ldots,w_{l(\rho)}$ with the partition $\rho$. For weights forming the partition $\rho$, the ordering of the vertices does not play a role. 

In a second step, we split up the original arrows $\alpha:p\to q$ by inserting a new vertex $pq$ and two new arrows $\alpha_1:p\to pq$ and $\alpha_2:pq\to q$ such that $w(\alpha)=w(\alpha_1)=w(\alpha_2)$. We write $Q'$ for these set of new vertices.

In a third step, we extend the ordering of the vertices $Q_0$ to the new vertices $V\cup W\cup Q'$ in such a way that again all arrows are pointing to larger vertices. Moreover, we want the ordering to satisfy $v_1>v_2>\ldots>v_{l(\lambda)}>v$ for every $v\in Q_0\cup W$.

\begin{defi} A weighted quiver $\tilde{\mathcal D}$ obtained in this way is called a $(\lambda,\rho)$-marked floor diagram of degree $d$ and genus $g$. 

The number of different $(\lambda,\rho)$-markings up to quiver automorphism which respect the ordering of $\mathcal D$ is denoted by $v_{\lambda,\rho}(\mathcal D)$. We define the multiplicity of $\tilde{\mathcal D}$ by
\[\mu_\rho(\tilde{\mathcal D})=\mu_\rho(\mathcal D)=\mu(\mathcal D)\prod_{i=1}^{l(\rho)}\rho_i.\]
 
\end{defi}
Note that $v_{\lambda,\rho}(\mathcal D)$ can be zero. Moreover, $\mu_\rho(\tilde{\mathcal D})$ only depends on the multiplicity of the original labeled floor diagram and on the partition $\rho$. 
\begin{thm}{\cite[Theorem 3.18]{fm}}
We have
\[N_{d,g}(\lambda,\rho)=\sum_{\mathcal D}\mu_\rho(\mathcal D)v_{\lambda,\rho}(\mathcal D)\]
where the sum is over all labeled floor diagrams of degree $d$ and genus $g$.

\end{thm}
Note that there is the same formula for the classical Gromov-Witten invariants $N_{d,g}$, see \cite[Theorem 1]{bf}. The special case $g=0$ is treated in \cite[Theorem 4.4]{abl}.

Let us consider the partition $\rho=\emptyset$ and $\lambda=(d)$. Let $\mathcal D$ be a labeled floor diagram of degree $d$ and genus $g$. Since we have $\mathrm{div}(v_{\max} )\leq 0$ for the largest vertex $v_{\max}$, it already follows that we only need to consider labeled floor diagrams such that $1-\mathrm{div}(v_{\max})=d$ for the largest vertex $v_{\max}$ and $1-\mathrm{div}(v)=0$ for all other vertices. Indeed, $V\cup W$ consists of only one vertex which is forced to be linked to the largest vertex. Note that the arrow which links the largest vertex to the new vertex does not affect the weight of the marked floor diagram.

If we additionally have $g=0$, this means that the labeled floor diagrams which we need to consider are certain weighted trees with root $v_{\max}$ such that $1-\mathrm{div}(v)=0$ for every vertex $v\neq v_{\max}$. We call these floor diagrams compatible with the partition $((d),\es)$. In this special case, it is possible to construct the combinatorial objects in a completely recursive way. As there is a similar description for the torus fixed points of the corresponding quiver moduli spaces, we review their construction.

Assume that $\mathcal D$ is such a diagram with largest vertex $\vm$ and assume that there are $k$ arrows of weights $a_k$ pointing to $\vm$ for $1\leq k\leq d-1$. Deleting the vertex $\vm$ and the corresponding $k$ arrows $\alpha_i$, we get $k$ diagrams $\cD_1,\ldots,\cD_k$.

Now the diagram $\cD_i$ is a labeled floor diagrams which is compatible with the partition $((a_i),\es)$. Indeed, if $i_{\max}$ is the vertex with $t(\alpha_i)=i_{\max}$, then we have $1-\dv(i_{\max})=a_i$ because we deleted the arrow $\alpha_i$ of weight $a_i$. For the other vertices of $\mathcal D_i$ we still have $1-\dv(v)=0$ because $\cD$ is a labeled floor diagram compatible with the partition $((d),\es)$. The other way around, given $k$ labeled floor diagrams $\cD_i$ compatible with $((a_i),\es)$, we get a labeled floor diagram $\cD$ compatible with $((\sum_{i=1}^{k}a_i+1),\es)$ by adding an extra vertex $\vm$ and $k$ extra arrows of weight $a_i$ from $i_{\max}$ to $\vm$. Then we have $1-\dv(i_{\max})=0$ and 
$$1-\dv(\vm)=1-(-1)\sum_{i=1}^ka_i=\sum_{i=1}^{k}a_i+1$$
where we understand $i_{\max}$ as a vertex of $\cD$.

Note that this construction is unique up to re-ordering $\{1,\ldots,k\}$. So what about the weights and markings? Since $\rho=\es$ and since we add $k$ arrows of weight $a_1,\ldots,a_k$, we have $$\mu_\es(\cD)=\prod_{i=1}^k\mu_\es(\cD_i)a_i^2.$$

To order the vertices of $\cD$ (which are $2d$ after splitting each arrow including the new arrows of weights $a_i$), we may split the set these vertices into $k$ parts of sizes $2a_1,\ldots,2a_k$. This just means that $2a_i$ vertices  correspond to the vertices of $\cD_i$. Thus, we only need to determine the number of orderings for each $\cD_i$ and obtain
\[v_{d,\es}(\cD)=\frac{(2d)!}{(2a_1)!\ldots(2a_k)!}\prod_{i=1}^kv_{a_i,\es}(\cD_i).\]
Taking the symmetries into account, for the corresponding Gromov-Witten-invariant we obtain the following recursive formulae.
\begin{thm}{\cite[Corollary 4.8, Theorem 4.9]{fm}}The following holds:
\begin{enumerate}
\item $$N_{d+1,0}((d+1),\es)=\sum_{k=1}^d\frac{(2d)!}{k!}\sum_{\substack{a_1+\ldots+a_k=d\\a_i>0}}\prod_{i=1}^k\frac{a_i^2N_{a_i,0}((a_i),\es)}{(2a_i)!}.$$
\item $N_{d,0}(\emptyset,(d))=d\cdot N_{d,0}((d),\emptyset)$.
\end{enumerate}
\end{thm}
Note that the second formula is immediate as the two invariants count the same diagrams, but the multiplicity differs by the factor $d$ for each diagram.

\begin{cor}\label{cor1} The numbers $N_{d,0}(\emptyset,(d))$ are determined by the recursion
$$N_{d+1,0}(\emptyset,(d+1))=(d+1)\cdot\sum_{k=1}^d\frac{(2d)!}{k!}\sum_{{a_1+\ldots+a_k=d}\atop{a_1,\ldots,a_k>0}}\prod_{i=1}^k\frac{a_i N_{a_i,0}(\emptyset,(a_i))}{(2a_i)!}.$$
\end{cor}

\subsection{MPS formula and localization for subspace quivers}
\noindent In order to determine the Euler characteristic of the moduli spaces $M_{d,2d+1}$, we apply the MPS-formula to the sink of the quiver $Q_m$ and combine it with the localization theorem, see Sections \ref{mpsform} and \ref{loc}.

We consider the linear form $\Theta\in (Q_m)_0$ defined by $\Theta_{i_k}=1$ for every $k=1,\ldots,m$ and $\Theta_j=0$ which is easily seen to be equivalent to the one fixed in Section \ref{point}. In particular, we have $$M_{d,2d+1}\simeq M_{\bf d}^{\Theta-{\rm st}}(Q_m)$$ for ${\bf d}=\sum_{k=1}^{2d+1}i_k+(2d+1)j$.

The MPS-formula for $Q_m$ applied to the unique sink can be expressed using the following quiver $\mathcal Q_m$. We define its vertices by
\begin{equation*}
(\mathcal Q_m)_0=\{i_k\mid k=1,\ldots,m\}\cup \{j_{(l,n)}\mid (l,n)\in\Nn^2\},
\end{equation*}
and the arrows by 
\begin{equation*}
(\mathcal Q_m)_1=\{\alpha_1,\ldots,\alpha_{l}:i_k\rightarrow j_{(l,n)},\forall\,l,n\in\Nn,\,k\in\{1,\ldots,m\}\}.
\end{equation*}
The level $l:(\mathcal Q_m)_0\to\Nn$ is given by
\begin{equation*}
l(j_{(l,n)})=l\,\,\text{for all } n\in\Nn,\,l(i_k)=1\,\,\text{for all } k\in\{1,\ldots,m\}
\end{equation*}
and induces the linear form $\Theta_l$ associating with every vertex its level. We consider the slope $\mu_l=\Theta_l/\kappa$ where $\kappa$ is defined by $\kappa(\bd)=\sum_{q\in (\mathcal Q_m)_0}l(q)\bd_q$ for $\bd\in \Nn(\mathcal Q_m)_0$. 

 Fix a dimension vector $\bd\in  \Nn(Q_m)_0$. Every (weighted) partition ${\bf P}=\sum^{t}_{l=1} lk_l$ of $\bd_j$, denoted by $\bp\vdash \bd_j$, defines a dimension vector $\bd({\bf P})\in\Nn(\mathcal Q_m)_0$ in the natural way, see Section \ref{mpsform}.

We write again $\mathcal M^{\Theta_l-\mathrm{st}}_{\bd(\bp)}(\mathcal Q_m)$ for the corresponding moduli space.
With this notation in place, the MPS-formula at the level of Euler characteristics applied to the sink $j$ can be expressed by
\begin{equation}\label{mps}
\chi(M_{\bf d}^{\Theta-{\rm st}}(Q_m)) = \sum_{\bp\vdash d_j}\chi(\mathcal M^{\Theta_l-\mathrm{st}}_{\bd(\bp)}(\mathcal Q_m))\prod_{l}\frac{(-1)^{k_l(l-1)}}{k_l!l^{2k_l}}.
\end{equation}
Thus we are left with the calculation of Euler characteristic of the moduli spaces $\mathcal M^{\Theta_l-\mathrm{st}}_{\bd(\bp)}(\mathcal Q_m)$. 

In the special case when the dimension vector $\bd$ satisfies $\bd_{i_k}=1$ for all $k=1,\ldots,m$, the dimension vector $d(\bp)$ is thin which means that this problem boils down to counting certain trees when applying the localization theorem of Section \ref{loc}.

Here, a spanning tree of the full subquiver $\mathrm{supp}(\bd(\bp))$ of $\mathcal Q_m$ is a subquiver of $\mathrm{supp}(\bd(\bp))$ which is a connected tree whose vertex set is $\mathrm{supp}(\bd(\bp))_0$. As alreday mentioned in Section \ref{loc}, every spanning tree defines a representation of $\mathrm{supp}(\bd(\bp))$ of dimension $\bd(\bp)$ in the obvious way. We call the spanning tree stable if the induced representation is stable with respect to $\Theta_l$. Note that there are $l$ arrows between the vertices $i_k$ and $j_{(l,n)}$ and every spanning tree contains at most one of them. By Theorem \ref{locforthin}, the Euler characteristic $\chi(\mathcal M^{\Theta_l-\mathrm{st}}_{\bd(\bp)}(\mathcal Q_m))$ is the number of stable spanning trees of $\mathrm{supp}(d(\bp))$. Recall that these are precisely the torus fixed points of $\mathcal M^{\Theta_l-\mathrm{st}}_{\bd(\bp)}(\mathcal Q_m)$ under the natural torus action.

Let us consider the case $\bd_{j}=d$, $m=2d+1$ and $\bd_{i_k}=1$ for all $k=1,\ldots,m$.  Since $(\ref{mps})$ does not yield a positive formula, on first sight it is not clear how to get a direct correspondence between the combinatoric objects, i.e. the stable spanning trees on the one side and labeled floor diagrams on the other side. But we do obtain a very easy description of the stable spanning trees (as well as we did for the labeled floor diagrams) which let us hope that there is a hidden direct correspondence between labeled floor diagrams and certain stable spanning trees. In turn, this makes the MPS-formula rather easy and it makes it possible to obtain an explicit formula for the Euler characteristic and thus for the Gromov-Witten invariants.

With a partition $\bp=\sum_llk_l$ we associate the full bipartite quiver $Q(\bp)$ with $\sum_{l\geq 1} k_l$ sinks and $2d+1$ sources where we have $k_l$ sinks of level $l$ and one (!) arrow from each source to each sink. Actually, $Q(\bp)$ is obtained from $\mathrm{supp} (\bd(\bp))$ by deleting $l-1$ arrows out of the $l$ arrows connecting a source to a sink of level $l$.
 
\begin{lem}\label{stabletree}
A spanning tree of $Q(\bp)$ is stable if and only if every sink of level $l$ has precisely $2l+1$ neighbours.
\end{lem}
This statement can be obtained when modifying the proof of \cite[Lemma 5.1]{wei2}, i.e. by an induction on the number of sinks. Thus we only give the idea. We say that a sink and a source are neighbors (or adjacent) if they are connected by an arrow. If the spanning tree has only one sink, the statement is clear. If it has more than one sink, it has a sink $j$, say of level $l$, which has precisely one adjacent source $i$ which has more than one neighbour. We can remove this sink together with all adjacent sources except $i$. By the stability condition, the number of these removed sources is forced to be $2l$. Again by the stability condition it can be seen that the remaining tree is a stable spanning tree of $Q(\bp-l)$ and we can apply the induction hypothesis which gives the claim.

We denote the number of stable spanning trees of $Q(\bp)$ by $n(\bp)$. Taking the possibilities of embeddings $Q(\bp)$ into $\mathrm{supp}(\bd(\bp))$ into account, we obtain:
\begin{thm}\label{coreuler}
We have
\begin{equation*}
\chi(M_{d,2d+1}) = \sum_{\bp\vdash d}n(\bp)\prod_{l\geq 1}\frac{(-1)^{k_l(l-1)}l^{k_l(2l-1)}}{k_l!}.
\end{equation*}
\end{thm}
We also obtain an explicit formula for the Euler characteristic
\begin{thm}
We have 
$$\chi(M_{d,2d+1})=\frac{(2d)!}{2d+1}\sum_{k=1}^n\frac{(-1)^{d-k}(2d+1)^k}{k!}\sum_{\substack{\sum_{i=1}^k a_i=d\\a_i>0}}\prod_{i=1}^k\frac{a_i^{2a_i-1}}{(2a_i)!}$$
\end{thm}

With each spanning tree, we associate a partition $k=\sum_l lk_l$ and a weight, which is induced by Theorem \ref{coreuler}, i.e.
\[\prod_{l\geq 1}\frac{(-1)^{k_l(l-1)}l^{k_l(2l-1)}}{k_l!}.\]
Now we can forget about the factor $\frac{1}{\prod_l k_l!}$ when we just count connected trees which have precisely $k_l$ sources of level $l$ and thus have $2l+1$ adjacent sinks. This means that we forget about the colouring of the sinks as each sink of a spanning tree corresponds to a sink of $Q(\bp)$. With each such tree we then associate a weight
\[\prod_{l\geq 1}(-1)^{k_l(l-1)}l^{k_l(2l-1)}.\]

Equivalently we can associate a weight $(-1)^{l-1}l^{2l-1}$ with each subquiver which has one sink of level $l$. Every source of a spanning tree also corresponds to a source of $\mathcal Q_m$. We forget about this correspondence for a while and let $y(x)$ be the generating function of such spanning trees taking weights into account. Now we can proceed analogously to \cite[Proposition 5.2]{wei2}. By Lemma \ref{stabletree}, every stable spanning tree can be obtained recursively. More detailed, we can fix an arbitrary stable spanning tree and a stable spanning tree which has only one sink of level $l$ and $2l+1$ sources and identify a source of the fixed tree with one of these $2l+1$ sources. We call a stable spanning tree with only one sink simple of level $l$.

If we do this algorithmically, we glue to a fixed source $b_1$ simple trees of level one, $b_2$ simple trees of level two and so on. 
Consider the power series
\[\Phi(x):=\exp\left(\sum_{l=1}^{\infty}\frac{x^{2l}}{(2l)!}(-1)^{l-1}l^{2l-1}\right).\]
Taking symmetries and the weights of each subquiver into account $y(x)$ satisfies the functional equation $y(x)=x \Phi(y(x))$. The number of such tree with $2d+1$ sources can now be obtained by the Lagrangian Inversion Theorem
\begin{eqnarray*}[x^{2d+1}]y(x)&=&\frac{1}{2d+1}[u^{2d}]\Phi(u)^{2d+1}\\&=&\frac{1}{2d+1}[u^{2d}]\sum_{n=0}^{\infty}\frac{1}{n!}\left((2d+1)\sum_{l=1}^{\infty}\frac{u^{2l}}{(2l)!}(-1)^{l-1}l^{2l-1}\right)^{n}\\&=&\frac{1}{2d+1}\sum_{k=1}^n\frac{(-1)^{d-k}(2d+1)^k}{k!}\sum_{\substack{\sum_{i=1}^k a_i=d\\a_i>0}}\prod_{i=1}^k\frac{a_i^{2a_i-1}}{(2a_i)!}.\end{eqnarray*}
Finally, we have to take into account that we  have to fix a source in the beginning in order to give this recursive description. For this choice, we have $2d+1$ possibilities. We also did not take the colourings of the sources into account, i.e. that every source corresponds to a source of the quiver $\mathcal Q_m$. In total this gives a factor $\frac{(2d+1)!}{2d+1}=(2d)!$. This shows the claim.

\begin{rem}
Thus we obtain a recursive description for the combinatorial objects which is very similar to the one for the label floor diagrams. Indeed, we can remove any source of a stable spanning tree which has more than one neighbour and obtain $k$ stable spanning trees corresponding to smaller dimension vectors of the same shape. One small difference between the two constructions is that this choice is not unique. A big one, which makes a comparison somehow difficult, is that the weights of the arrows of the labeled floor diagrams are ascending in the direction of the largest vertex.
\end{rem}


\begin{thebibliography}{99}

\bibitem{abl} Arroyo, A., Brugallé, E., De Medrano, L. L.: Recursive formulas for Welschinger invariants of the projective plane. International Mathematics Research Notices \textbf{2011}(5), 1107--1134 (2010).
\bibitem{bgp}
Bernstein, I.~N., Gelfand, I.~M., Ponomarev, V.~A.: Coxeter functors, and Gabriel's theorem. Russian Mathematical Surveys \textbf{28}(2), 17--32 (1973).
\bibitem{BG} Block, F.,  G\"ottsche, L.: Fock spaces and refined Severi degrees.
International Mathematics Research Notices \textbf{2016}(21), 6553--6580 (2015).

\bibitem{bf} Brugallé, E., Mikhalkin, G.: Enumeration of curves via floor diagrams. Comptes Rendus Mathematique \textbf{345}(6), 329--334 (2007).
\bibitem{ER} Engel, J., Reineke, M.: Smooth models of quiver moduli. Mathematische Zeit\-schrift \textbf{262}(4), 817--848 (2009).
\bibitem{FR} Franzen, H., Reineke, M.: Cohomology Rings of Moduli of Point Configurations on the Projective Line. Preprint, arxiv:1611.01092 (2016).
\bibitem{fm} Fomin, S., Mikhalkin, G.: Labeled floor diagrams for plane curves. Journal of the European Mathematical Society \textbf{12}, 1453--1496 (2010). 
\bibitem{gps} Gross, M., Pandharipande, R., Siebert, B.: The tropical vertex. Duke Mathematical Journal \textbf{153}(2), 297--362 (2010).
\bibitem{HMSV} Howard, B., Millson, J., Snowden, A., Vakil, R.: The equations for the moduli space of n points on the line.  Duke Mathematical Journal \textbf{146}(2), 175--226 (2009).
\bibitem{HMSV1} Howard, B.,  Millson, J., Snowden, A.,  Vakil, R.: A description of the outer automorphism of S6, and the invariants of six points in projective space.
 Journal of Combinatorial Theory, Series A, \textbf{7}(115), 1296--1303 (2008).
\bibitem{kin}
King, A.: Moduli of representations of finite-dimensional algebras. Quarterly Journal of Mathematics \textbf{45}(4), 515-530 (1994).
\bibitem{KM} Kontsevich, M., Manin, Y.:  Gromov-Witten classes, quantum cohomology, and enumerative geometry. Communications in Mathematical Physics \textbf{164}(3), 525--562 (1994).
\bibitem{KS} Kontsevich, M.,  Soibelman, Y.: Cohomological Hall algebra, exponential Hodge structures and motivic Donaldson-Thomas invariants.  Communications in Number Theory and Physics  \textbf{5}, 231--352 (2010).
\bibitem{mps} Manschot, J., Pioline, B., Sen, A.: Wall crossing from Boltzmann black hole halos. Journal of High Energy Physics \textbf{7}, 1--73 (2011).
\bibitem{MeR} Meinhardt, S., Reineke, M.: Donaldson-Thomas invariants versus intersection cohomology of quiver moduli. Preprint, arXiv:1411.4062 (2014). 
\bibitem{MFK} Mumford, D.; Fogarty, J.; Kirwan, F.:
Geometric invariant theory. Third edition. Ergebnisse der Mathematik und ihrer Grenzgebiete (2), 34. Springer-Verlag, Berlin, 1994.
\bibitem{rsw} Reineke, M., Stoppa, J., Weist, T.: MPS degeneration formula for quiver moduli and refined GW/Kronecker correspondence. Geometry \& Topology \textbf{16}(4), 2097--2134 (2012).
\bibitem{rw} Reineke, M., Weist, T.: Refined GW/Kronecker correspondence. Mathematische Annalen \textbf{355}(1), 17--56 (2013).
\bibitem{HNS} Reineke, M.: The Harder-Narasimhan system in quantum groups and cohomology of quiver moduli. Inventiones Mathematicae \textbf{152}(2), 349--368 (2003).
\bibitem{RModuli} Reineke, M.: Moduli of representations of quivers.
 Trends in representation theory of algebras and related topics, EMS Series of Congress Reports, European Mathematical Society, Zürich, 2008, 589--637. 
\bibitem{RSm} Reineke, M.: Quiver moduli and small desingularizations of some GIT quotients. Representation theory -- current trends and perspectives, EMS Series of Congress Reports, European Mathematical Society, Z\"urich, 2017, 613--635.
\bibitem{RS} Reineke, M.,  Schr\"oer, S.: Brauer groups for quiver moduli.
 Algebraic Geometry \textbf{4}(4), 452--471 (2017).

\bibitem{wei}  Weist, T.: Localization of quiver moduli spaces. Representation Theory \textbf{17}(13), 382-425 (2013).
\bibitem{wei2}Weist, T.: On the Euler characteristic of Kronecker moduli spaces. Journal of Algebraic Combinatorics \textbf{38}(3), 567-583 (2013).
\end{thebibliography}
\end{document}